\documentclass{amsart}
\usepackage{psfig,latexsym,amssymb}
\def\proof{\medskip\noindent{\sc Proof. }}

\def\dwreath{\Delta(B_{2n})/S_2\wr S_n }

\def\EOP{\hfill$\Box$}
\def\natnum{\hbox{\rm I\kern-.17em N}}
\def\integ{\hbox{\rm Z\kern-.3em Z}}
\def\real{R}
\def\dwr{ \Delta (B_{lm})/S_l\wr S_m }
\def\dring{ k[x_1,\dots ,x_{lm}]^{S_l\wr S_m} }

\newtheorem{thm}{Theorem}[section]
\newtheorem{lem}{Lemma}[section]
\newtheorem{prop}{Proposition}[section]
\newtheorem{defn}{Definition}[section]

\newtheorem{qn}{Question}[section]
\newtheorem{examp}{Example}[section]

\newtheorem{rk}{Remark}[section]

\begin{document}

\title[Properties of the quotient complex $\Delta (B_{lm})/S_l \wr
S_m$]{A partitioning and related properties for
the quotient complex $\Delta(B_{lm})/S_l\wr S_m$}
\author{Patricia Hersh}
\address{Department of Mathematics,
University of Michigan, Ann Arbor, Michigan 48109-1109}
\curraddr{
        Department of Mathematics, 
        University of Michigan,
	525 East University Ave.
	Ann Arbor, Michigan 48109-1109}
\email{plhersh@umich.edu}
\thanks{The author was supported by an NSF postdoctoral research
fellowship.}

\begin{abstract}
We study the quotient complex $\dwr $ as a means of deducing facts 
about the ring $\dring $.  It is shown in [He] that $\dwr $ is 
shellable when $l=2$, implying Cohen-Macaulayness of 
$k[x_1,\dots ,x_{2m}]^{S_2\wr S_m}$ for any field $k$.
We now confirm for all pairs $(l,m)$ with $l>2$ and $m>1$ that 
$\dwr $ is not Cohen-Macaulay over $\integ /2\integ $,
but it is Cohen-Macaulay over fields of characteristic 
$p>m$ (independent of $l$).  This 
yields corresponding characteristic-dependent
results for $\dring $.  We also prove that $\dwr $ and the links of many
of its faces are collapsible, and we give a partitioning for
$\dwr $.  
\end{abstract}

\maketitle

\section{Introduction}\label{introsec}

Let $B_n$ denote the Boolean algebra of subsets of $\{ 1,\cdots ,n\} $ 
ordered by inclusion. The natural symmetric 
group action on $\{ 1,\cdots ,n \} $ induces a rank-preserving, 
order-preserving action on $B_n$.
Likewise, the wreath product of symmetric groups $S_l \wr S_m \subset S_{lm}$ 
acts on the Boolean algebra $B_{lm}$.  (Recall that $S_l\wr S_m$ is
the subgroup of $S_{lm}$ of order $(l!)^m m!$ which permutes the values
$il+1,\dots ,(i+1)l$ among themselves for each $0\le i< m$ and also 
wholesale permutes these $m$ sets of size $l$.)
This induces an $S_l\wr S_m$-action 
on chains $\hat{0} < u_0 \cdots < u_i < \hat{1}$ of comparable poset 
elements, i.e. on faces in the order complex $\Delta(B_{lm}) $.
The action on chains 
gives rise to a quotient cell complex, denoted $\Delta(B_{lm})/S_l\wr
S_m$, which consists of the $S_l\wr S_m$-orbits of the order complex faces.  
As a word of caution, the quotient complex $\Delta(B_{lm})/S_l\wr Sm $ does
not coincide with the order complex of the quotient poset $B_{lm}/S_l\wr
S_m$ (cf. [BK] for a study of which quotient complexes are order complexes
of quotient complexes), because there are covering relations 
$u\le v, u'\le v'$ in $B_{lm}$ belonging to 
distinct orbits despite having $u'=gu$ and $v'=g'v$ for some $g,g'\in S_l\wr
S_m$.

We will rely on results of Stanley, Hochster-Eagon, Reiner, Bj\"orner and
Garsia-Stanton to transfer properties of the quotient complex
$\Delta(B_{lm})/S_l\wr S_m$ into algebraic facts about the
subring of invariant polynomials $\dring $.  Section 2 will review these
results about subrings of invariant polynomials, quotient complexes and 
more generally about simplicial posets from [Bj], [GS], [HE], [Re] 
and [St3].
Sections 3 and 4 follow up on previous work in [He], where a lexicographic
shelling was given for $\Delta (B_{2m})/S_2\wr S_m$.  In section 3, we 
show that $\dwr $ is not Cohen-Macaulay over the integers mod 2 whenever 
$l>2$ and $m>1$, by exhibiting local 2-torsion.  (The
situation is trivial whenever $l=1$ or $m=1$.)
Section 4 shows that $\Delta (B_{lm})/S_l \wr 
S_m $ and many of its links are collapsible, and finally we provide a 
partitioning for $\dwr $ in Section 5.

One theme that runs
throughout this paper is the use of ideas typically associated 
(at least implicitly) to 
lexicographic shellings to deduce properties related to shellability
for complexes that are not shellable; in particular, we give
collapsibility, Cohen-Macaulayness (for certain field characteristics) and 
partitionability results for $\dwr $.  In theory, our partitioning for
$\dwr $ gives a Hilbert series expression for $k[x_1,\dots ,
x_{lm}]^{S_l\wr S_m}$, but it would be desirable to find a simpler
expression.  Our partitioning for $\dwr $ is very similar to the
latter half of the (very complicated) partitioning argument used in [He] for 
$\Delta (\Pi_n)/S_n$; one of our goals was to simplify that argument.  

It remains open for $l>2$ to determine for which field characteristics $p$
such that $2<p \le m$ the ring $\dring $ is Cohen-Macaulay.  Garsia-Stanton
showed in [GS] how to deduce Cohen-Macaulayness over fields of 
characteristic $p$ from partitionings in which
$p$ does not divide the determinant of the incidence matrix.  We hope
that our work may help with the resolution of this question.

\section{Simplicial posets, quotient complexes and subrings of invariant
polynomials}

Boolean cell complexes were defined as follows in [Bj] and [GS]:

\begin{defn}
A regular cell complex is {\bf boolean} if every lower-interval
in its face poset is a Boolean algebra, namely if each cell has the 
combinatorial type of a simplex.  
\end{defn}

\medskip
Stanley studied their face posets, which he called simplicial
posets, in [St3].  People often use the term simplicial poset to
mean either the face poset or the cell complex itself; we will  
reserve the term simplicial poset exclusively for the face posets, to 
emphasize the distinction between a boolean cell complex and the 
order complex of its face poset.

One may think of the cells in a boolean cell complex as simplices,
but unlike in simplicial complexes, multiple faces may
have the same set of vertices.  As a result,
two faces may overlap in a simplicial complex rather than simply in
a face.  We refer to $i$-cells as $i$-faces,
0-cells as vertices, and call cells of top dimension facets.  Our
interest is in a particular class of boolean cell complexes, namely 
the quotient complexes $\Delta/G $ made up 
of the $G$-orbits of faces
in a simplicial complex $\Delta $ when a group $G$ acts simplicially 
on the faces of $\Delta $.

Stanley defined the face ring $k[P]$ for a simplicial poset $P$ in [St3]
by taking the faces in a boolean cell complex (or equivalently the elements
in its face poset $P$) as the generators of a polynomial ring over a field
$k$ and giving the generators 
the following three types of relations: 
\begin{enumerate}
\item
$xy$ if there is no face containing both $x$ and $y$
\item
$xy - (x\wedge y)\left(\sum_{z \in lub(x,y)} z\right)$, where $lub(x,y)$ 
denotes the set of least upper bounds of $x$ and $y$
\item
$\hat{0} - 1$
\end{enumerate}

Stanley proved the following in 
[St3], using facts about algebras with straightening laws.

\begin{thm}[Stanley] 
The face ring $k[P]$ of a Cohen-Macaulay simplicial poset $P$ is
a Cohen-Macaulay ring.
\end{thm}

Let us denote the face ring of the face poset of a 
quotient complex $\Delta /G$ by 
$k[\Delta /G]$.  In [Re], Reiner established
the following connection between face rings
of quotient complexes and subrings of invariant polynomials (cf. [St4, p.
53] for the definition of $k[\Delta ]$, or specialize the above definition
to simplicial complexes).

\begin{thm}[Reiner] The rings $k[\Delta /G] $ and $k[\Delta ]^G$ are 
isomorphic.
\end{thm}

Reiner also showed (unpublished) that Cohen-Macaulayness
for subrings of invariant polynomials for 
face rings of certain quotients of type A Coxeter complexes transfers
to Cohen-Macaulayness of other subrings of invariant polynomials.  A proof
of the following result has been provided by 
Reiner in an appendix.

\begin{thm}[Reiner]
\label{appendix-result}
If $G\subset S_n$ and $k[\Delta (B_n)]^G $ is 
Cohen-Macaulay over a field $k$, then $k[x_1,\dots ,x_n]^G$ is  
Cohen-Macaulay over the same field $k$. 
\end{thm}

In [Bj4], Bj\"orner established a notion 
of shellability for boolean cell complexes (stated slightly differently than
below) and noted that it implies Cohen-Macaulayness.

\begin{defn}[Bj\"orner]
A pure boolean cell complex is {\bf shellable} if the 
facets may be ordered $F_1,\dots ,F_k$ so that 
$F_j \cap (\cup_{i=1}^j F_i )$ is pure of codimension one for each
$1<j\le k$. 
\end{defn}

\medskip
Just as in the case of simplicial complexes, this 
is equivalent to requiring there to be a unique minimal new face at
each facet insertion.

\begin{prop}[Bj\"orner]\label{goofcm}
If a pure boolean cell complex is shellable, then the underlying
topological space is Cohen-Macaulay (over any field).
\end{prop}

\medskip
We will also use the following result of Hochster and 
Eagon to get at the Cohen-Macaulayness of $k[x_1,\dots ,x_{lm}]^{S_l\wr 
S_m}$ for relatively large field characteristics.

\begin{thm}[Hochster-Eagon] If $\Delta $ is a Cohen-Macaulay simplicial
complex and the characteristic of $k$ 
does not divide $|G|$, then $k[\Delta ]^G$ is a Cohen-Macaulay ring.
\end{thm}

In discussing which complexes $\Delta (B_{lm})/S_l\wr S_m$ are shellable,
we will make use of the fact that 
$\Delta (B_{lm})/S_l\wr S_m$ is balanced.  Recall that a 
boolean cell complex of dimension $d-1$ 
is {\bf balanced } if there is a map
$\kappa : V(\Delta ) \rightarrow \{ 1,\dots ,d\} $ that colors the vertices
with $d$ colors so that no two vertices in the same face are the 
same color.  We refer to the set of colors for the vertices in a face 
as the {\bf support} of the face.  Notice that the order complex of a 
finite, graded poset is balanced by poset rank.
One nice feature of balanced complexes
is that their face rings have very explicit linear systems of 
parameters (l.s.o.p.'s), namely the face ring of a
balanced $(d-1)$-dimensional complex $\Delta $ has linear system of parameters
$\theta_1,\dots ,\theta_d$ in which $\theta_i = \sum_{v: \kappa (v)=i} v$
(cf. [St4]).  

If a complex $\Delta $ of dimension $d-1$ is shellable and $k[\Delta ]$ has 
linear system of parameters $\theta_1,\dots ,\theta_d$, then $k[\Delta ] 
= \coprod_{\nu \in X} \nu \cdot k[\theta_1,\dots ,\theta_d]$ and the set
$X $ of minimal faces in the shelling is 
a $k$-basis for $k[\Delta ]/(\theta_1,\dots ,\theta_d)$ (cf. [St4]).  In 
this case, $X$ is called a {\bf basic set} for $k[\Delta]$.  
Garsia and Stanton use shellings and certain types of partitionings as a 
means for constructing basic sets for 
rings $k[\Delta /G]$ and for related subrings of 
invariant polynomials in [GS];  we
follow their notation in the remainder of this section.

If $c$ is a 
face of $\Delta $ consisting of vertices $x_{i_1},\dots ,x_{i_r}$, then 
denote by $x(c)$ the monomial $x_{i_1}\cdots x_{i_r}$ in the 
face ring $k[\Delta ]$.  When a group $H$ acts on $\Delta $, the Reynold's
operator $R^H$ acts on $k[\Delta ]$ by
$$R^H (x(c)) = \frac{1}{|H|}\sum_{h\in H} hx(c) = \frac{1}{|H|}\sum_{h\in H}
x(hc).$$  A set of chain monomials $\{ x(b) | b\in B\} $ given by a 
collection $B$ of chains in a poset $P$ is called a {\bf basic set} if every
element $Q$ of the Stanley-Reisner ring $k[\Delta ]$ has a unique expression
$$Q = \sum_{b\in B} x(b) Q_b (\theta_1,\dots ,\theta_d)$$ where the 
coefficients $Q_b(\theta_1,\dots ,\theta_d )$ are polynomials with rational
coefficients in the variables $\theta_1,\cdots ,\theta_d $.  This yields a 
Hilbert series expression $$Hilb (k[\Delta ],\lambda ) = \left( \prod_{i=1}^d
\frac{1}{1-\lambda^{\deg (\theta_i )}}\right) \left(\sum_{b\in B}\lambda^{
\deg (x(b))}\right).$$  All Cohen-Macaulay posets have such
basic sets.

\begin{thm}[Garsia-Stanton]\label{gs}
If $\Delta /H$ has a shelling $F_1,\dots ,F_k $ where $G_j $ is the unique
minimal new face in $F_j \setminus (\cup_{i<j} F_i )$ and $b_j$ 
is a representative of the orbit $G_j$ within $\Delta $,
then the orbit polynomials $R^H x(b_i) $ form a basic set for the 
subring of invariant polynomials $k[\Delta ]^H$, implying Cohen-Macaulayness
over any field.  
\end{thm}

\medskip
When a subgroup $G$ of the symmetric group $S_n$ acts on
the boolean algebra $B_n$ in a rank-preserving, order-preserving fashion,
then Garsia and Stanton proved in [GS] that basic sets for
$k[\Delta (B_n)]^G$ transfer to basic sets for $k[x_1,\dots ,x_n]^G $ 
and that certain types of partitionings (including all shellings) give rise to
basic sets.  We state their result in Theorem ~\ref{gsshell}, 
but first we give a definition it will use.

\begin{defn}
The {\bf incidence matrix} of a partitioning is a matrix with rows indexed
by facets and columns indexed by the minimal faces in the partitioning.
If $G_j\subseteq F_i$ then $A_{i,j}=1$ and otherwise $A_{i,j}=0$.
\end{defn}

\medskip
The incidence matrix for a partitioning coming from a shelling 
is always upper triangular with 1's on the diagonal, hence nonsingular
(over any field).  Other partitionings may yield incidence matrices that are 
singular over finite fields of sufficiently small characteristic.  It is 
possible to
construct partitionings with singular incidence matrices for 
Cohen-Macaulay complexes (personal communication of Reiner), so one cannot
conclude non-Cohen-Macaulayness by obtaining a singular incidence matrix.

\begin{thm}[Garsia-Stanton]\label{gsshell}
Let $G\subset S_n$ act as above
and let $[G_1,F_1]\cup\cdots\cup [G_k,F_k]$ be a partitioning for
$\Delta (B_n)/G$ with nonsingular incidence
matrix.  Then $x(G_1),\dots ,x(G_k)$ form a basic set for
$k[\Delta (B_n)/G]$, w.r.t. \!the l.s.o.p. \!$\theta_1,\dots ,
\theta_{n-1}$ given by the balancing.
Sending $\theta_i$ to the elementary symmetric function $e_i$ and 
$G_j = S_1 \subset S_2 \subset \cdots \subset S_r$ to the 
product $x_{S_1}x_{S_2}\cdots x_{S_r}$, in which $x_S = \prod_{i\in S}x_i$
yields a basic set for $k[x_1,\dots ,x_n]^G$.
\end{thm}


\section{Shellability and Cohen-Macaulayness results}

Using a lexicographic shellability criterion for pure, balanced complexes,
it is shown in [He] that $\Delta(B_{2n})/S_2\wr S_n $ is shellable.
Below we will describe the lexicographic order that led to this shelling,
but we refer readers to [He] for the proof that it does indeed give a 
shelling.

The following chain-labeling 
for $\Delta(B_{2m})/S_2\wr S_m $ gives a lexicographic shelling:
label the covering relation
$\{ \sigma_1,\dots ,\sigma_{i-1}\} \prec \{ \sigma_1,\dots ,\sigma_{i-1},
\sigma_i \}$ in the poset $B_{2m}$ 
with the label $\sigma_i\in \{ 1,\dots ,2m\} $, 
recording the insertion of $\sigma_i$.  Thus, the
saturated chain $\emptyset\prec\{\sigma_1\}\prec
\cdots\prec \{\sigma_1,\dots ,\sigma_{2m}\} $ is labeled $\sigma_1\cdots
\sigma_{2m}\in S_{2m}$.
The facets in $\Delta(B_{2m})/S_2\wr S_m$ are the orbits of the
saturated chains in $B_{2m}$, and by convention we label 
each of these orbits with
lexicographically smallest permutation among the labels for
members of the orbit.
This chain-labeling gives a $CC$-shelling, in the 
sense developed for posets by Kozlov in [Ko1] and extended to pure, balanced 
complexes in [He]. (Hultman recently further generalized the lexicographic
shellability criterion of [He] to non-pure balanced complexes in [Hu].)

The labels for the orbit representatives turn out to be
the permutations of $1,\dots ,2m$ which do not have any inversion pairs
$(2i-1,2i)$ or $(2i-1,2i+1)$, namely permutations in which
the odd numbers 
appear in increasing order and each odd number comes earlier than its
even successor.

\begin{examp}
\rm
The orbit representatives for $\Delta(B_6)/S_2\wr S_3 $, listed in
lexicographic order, are 
$123456$,
$1235\bullet 46$,
$123\circ 56\bullet 4$,
$13\bullet 2456$,
$13\bullet 25\bullet 46$,
$13\bullet 256\bullet 4$,
$1\circ 34\bullet 256$,
$1\circ 345\bullet 26$,
$1\circ 3456\bullet 2$,
$135\bullet 246$,
$13\circ 5\bullet 26\bullet 4$,
$1\circ 35\bullet 4\bullet 26$,
$1\circ 35\bullet 46\bullet 2$,
$13\circ 56\bullet 24$, and 
$1\circ 3\circ 56\bullet 4\bullet 2$.
Hollow dots denote ascents which behave topologically like descents
and filled-in dots indicate traditional descents.  The minimal new 
face for a facet is the union of the ranks of the hollow dots and 
the ranks of the filled-in dots.  For instance,
the swap ascent in $1\circ 3456\bullet 2$ comes from a
codimension one face skipping rank 1 in the intersection of
$134562$  with $132564$, resulting from
the fact that $312564$ is in the same orbit as $134562$.
\end{examp}

\begin{figure}[h]
\begin{picture}(150,62)(-65,0)
\psfig{figure=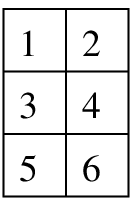,height=2cm,width=1.3cm}
\end{picture}
\caption{Labeled boxes acted upon by $S_2\wr S_3$}
\label{wreathbox}
\end{figure}

\medskip
To describe the group $S_2\wr S_m$ (and more generally $S_l\wr S_m$), 
let us first place the numbers $1,\dots ,2m$ (resp. 
$1,\dots ,lm $) in a $2\times m$ (resp. $l\times m$) table, by
sequentially inserting the numbers from left to right in each row, proceeding
from one row to the next from top to bottom, as in Figure ~\ref{wreathbox}. 
The elements of $S_2\wr S_m$ (resp. $S_l\wr S_m$) may then be described as 
the permutations in $S_{2m}$ (resp. $S_l\wr S_m$) which permute the 
numbers within each row and then permute the set of rows.
More formally, each element of $S_2\wr S_m$ is a 
composition of $\pi_2\circ \pi_1\in S_{2m}$ 
in which $\pi_1 = (12)^{e_1}
(34)^{e_2}\cdots (2m-1,2m)^{e_m}$ for some vector $(e_1,\dots ,e_m)\in 
\{0,1\}^m$ and where $\pi_2$ is obtained from some
$\pi\in S_m$ by requiring
$\pi_2 (2i) = 2\pi (i)$ and $\pi_2(2i-1) = 2\pi(i)-1$ for $1\le i\le m$
(and of course there is a similar definition for $l>2$).

The lexicographic shelling for $\Delta (B_{2m})/S_2\wr S_m $ may be 
combined with results of Stanley and Reiner (recalled in Theorems 2.1 and 
2.2, respectively) to obtain
Cohen-Macaulayness for $k[x_1,\dots ,x_{2m}]^{S_2\wr S_m}$ independent of
field characteristic (or equivalently over the integers), as is 
noted in [He].  
When $char(k)=0$, this is a special case 
of a result from [HE], but the shelling also allows coefficients in fields
of finite characteristic or the integers.
By Theorem 2.6, the
lexicographic shelling for $\dwreath $ also yields a basic set for the
subring $k[x_1,\dots ,x_{2m}]^{S_2\wr S_m}$ of invariant polynomials.
A simple description of which descent sets
occur in the lexicographic shelling would be desirable in that it would
yield a nice description of 
these basic sets (and in turn a nice Hilbert series
expression).

The story is more subtle for 
$\Delta(B_{lm})/S_l \wr S_m$ when $l$ is greater than $2$.  It is 
observed in [He] that these complexes cannot be shellable when 
$l>2$, by a Molien series computation which shows that the Hilbert 
series disagrees with the expression that would result from applying
Theorem 2.5 (recalled from [GS]) to any potential shelling order.  Now 
we construct
explicit faces whose links have 2-torsion and give partial results 
regarding the question of for which coefficient fields is $\dwr $ 
Cohen-Macaulay.  In particular, the fact that $\Delta (B_{lm})$ is a
triangulation of a sphere immediately implies (via a result 
of Hochster and Eagon [HE]) that $\dwr $is Cohen-Macaulay for coefficient 
fields of characteristic $p$ so long as $p$ does not divide 
$|S_l\wr S_m|$, i.e. for primes $p$ larger than $\max (l,m)$.  We will 
do slightly better, showing
Cohen-Macaulayness for $p>m$, regardless of how large $l$ grows.
We also show that $\dwr $ and the links of many faces are collapsible,
restricting how local $p$-torsion in lower homology might arise.

For each pair $(l,m)$ with $l>2$ and $m>1$, we will 
provide a face $F$ such that $\rm{lk} (F)$ has dimension two and also has
homology group $H_1(\Delta ,\integ ) = \integ /2\integ $, precluding 
Cohen-Macaulayness.  First
consider the link of the face $F = \emptyset\subseteq \{ 1,4\} \subseteq
\{ 1,2,4,5\} \subseteq \{ 1,2,3,4,5,6\} $ in 
$\Delta(B_6)/S_3\wr S_2$.  
Notice that $\rm{lk}(F)$ has 3 vertices, 6 edges and 4
2-simplices, and that the underlying topological
space is the real projective plane $\real P_2$.  

\begin{prop}\label{non-cm}
\rm
The quotient complex $\Delta (B_{lm})/S_l\wr S_m$ is not Cohen-Macaulay over
$\integ /2\integ $ for $l>2$ and $m>1$, hence $\integ /2\integ [x_1,\dots ,
x_{lm}]^{S_l\wr S_m}$ is not Cohen-Macaulay for such pairs $(l,m)$.
\end{prop}

\proof
One gets $\real P_{l-1}$ as the link of a face in
$\Delta (B_{2l})/S_l\wr S_2$, as follows: let us call all the letters
in the first ``row'' 1 and all the letters in the second ``row'' 2 
(since the letters in a row are all interchangeable) and then take the face 
$$F=\emptyset \subseteq \{ 1, 2\} \subseteq \{ 1^2 ,2^2\}\subseteq\cdots
\subseteq \{ 1^l ,2^l \} .$$  Note that the link of this face in 
$\Delta (B_{2l})/S_l \times S_2 $ is a sphere, because $\Delta(B_{2l})/
S_l\times S_2$ is lexicographically shellable (as shown by Garsia and 
Stanton in [GS]), and the restriction of this shelling to $\rm{lk}(F)$ in
$\Delta (B_{2l})/S_l \times S_2 $ has one decreasing chain.  We obtain
the desired link in $\Delta (B_{2l})/S_l\wr S_2$ 
by gluing together pairs of antipodal faces in this sphere (i.e. by 
identifying faces in which the two classes of objects are exchanged);
thus we obtain projective space in a completely natural fashion.  This
link also sits inside $\Delta (B_{lm})/S_l\wr S_m$ as the link of a 
larger face.

The conclusion about $\integ /2\integ [x_1,\dots ,
x_{lm}]^{S_l\wr S_m}$ for such pairs $(l,m)$ follows from the same 
reasoning used for other coefficient fields in Proposition ~\ref{cmfield}.
\EOP

\begin{prop}\label{cmfield}
The quotient complex $\dwr $ and consequently the ring $\dring $ 
is Cohen-Macaulay over fields of 
characteristic $p$ whenever $p>m$.
\end{prop}
\proof
It is shown that $\Delta (B_{lm})/S_l \times \cdots \times S_l $ 
is shellable (and hence Cohen-Macaulay over any field) in [GS].
Note that $\Delta (B_{lm})/S_l \wr S_m $ is the quotient of a Cohen-Macaulay
complex by an $S_m$-action, by virtue of the isomorphism
$\Delta (B_{lm})/S_l \wr S_m \cong
\left( \Delta (B_{lm})/S_l\times  \cdots \times S_l\right) /S_m$.
Thus, one may apply the result of Hochster and Eagon [HE], recalled in
Theorem 2.4, to conclude that there is no $p$-torsion unless 
$p$ divides $|S_m|$, i.e. unless $p\le m$.  

Stanley's result from [St3] that face rings of Cohen-Macaulay simplicial
posets are Cohen-Macaulay then tells 
us that $k[\Delta (B_{lm})/S_l\wr S_m ]$ is
Cohen-Macaulay for $p=char(k) > m$, but $k[\Delta (B_{lm})/S_l\wr S_m ]
\cong k[\Delta (B_{lm})]^{S_l\wr S_m} $ by Theorem 2.2.
Now we apply Theorem 2.3 to conclude
that $k[x_1,\dots ,x_{lm}]^{S_l \wr S_m} $ is also Cohen-Macaulay for
$p=char(k)>m$.
\EOP

\begin{qn}
Is there any local $p$-torsion in lower homology for $2 < p \le m$?  
\end{qn}
\medskip
In Section 5, we will give a partitioning for $\Delta (B_{lm})/S_l\wr S_m$,
and if the determinant of the incidence matrix for this partitioning were
not divisible by a prime $p$, then one could conclude 
Cohen-Macaulayness of $\dring $ for $char(k)=p$.  We 
suspect that $p$-torsion for primes larger than 2 would already appear in
$\rm{lk}(\emptyset < \{ 1,2,3\} < \{ 1^2,2^2,3^2 \} < \{ 1^3,2^3,3^3 \} )$
if it ever occurs.  The determinant of the incidence matrix for our 
partitioning of this link is $2^3\cdot 3^5$, strongly suggesting (but not 
confirming) there is local 3-torsion present.

\begin{rk}
The directed graph complexes studied by Kozlov in
[Ko2] have faces whose links are isomorphic to $\rm{lk}(\emptyset < 
\{ 1,\dots ,m\} < \cdots < \{ 1^l,\dots ,m^l \} )$ in 
$\Delta(B_{lm})/S_l\wr S_m$ for any pair $(l,m)$, and 
hence there is local $2$-torsion arising just as in Proposition 3.1.
Kozlov previously determined by computer
that the directed graph complexes have local $2$-torsion.
\end{rk}

\section{Collapsibility of $\Delta (B_{lm})/S_l\wr S_m$ and links of 
many faces}

This section proves that $\Delta (B_{lm})/S_l\wr S_m$ and 
the links of many of its faces are collapsible.  The discussion of links
is included in the hope that this might
shed some light on the question of when the complexes are Cohen-Macaulay
(i.e. for which field characteristics $p$ such
that $2<p\le m$).
The collapsibility proofs are a relaxation of the sort of 
argument typically used to produce lexicographic shellings in that we
will show that the intersection $F_j \cap (\cup_{i<j} F_i )$ of each 
facet $F_j$ with the union of earlier ones is collapsible by exhibiting
a topological ascent in each $F_j$, yielding a cone point in each 
intersection $F_j \cap (\cup_{i<j} F_i)$.  

\begin{thm}
The quotient complex $\dwr $ is collapsible.
\end{thm}

\proof
Let us first order the saturated chains in $B_{lm}$ lexicographically, just as
in the lexicographic shelling for $\Delta (B_{2m})/S_2\wr S_m$, 
and then choose the 
lexicographically earliest saturated chain in each $S_l\wr S_m$-orbit as the
orbit representative.  Now we build up the quotient
complex by sequentially inserting
facets of $\dwr $ in the resulting lexicographic order $F_1, \dots , F_r$.
We will prove collapsibility by showing 
that each intersection $F_j\cap (\cup_{i<j} F_i)$ for $j>1$
has a cone point so 
that collapsibility is preserved with each facet insertion as we sequentially
build the complex, since clearly $F_1$ is itself collapsible.  

Let us encode the permutations in $S_{lm}$ which label the saturated chain 
orbit representatives as words of 
content $\{ 1^l,2^l,\dots ,m^l \}$ by replacing the label
$rl+s$ with the label $r+1$ for each $0\le r < m, 0<s<l$.  
Notice that this map is a bijection between
permutations in $S_{lm}$ which are lexicographically smallest in their
$S_l\wr S_m$-orbit and words of content $\{ 1^l,2^l,\dots ,m^l \} $ in which 
the first appearance of $i$ precedes the first appearance of $j$ for each 
pair $1\le i<j\le m$.  Note that any descent in 
the labels on a saturated chain orbit $F_j$
may be replaced by a lexicographically smaller ascent to get the label
for a lexicographically earlier saturated chain orbit $F_i$ such 
that $F_i$ and $F_j$ share
a codimension one face obtained by omitting the descent from $F_j$.  
We will show that $F_j\cap (\cup_{i<j} F_i)$ has a cone point at
the element $u\in F_j$ just preceding the final appearance of 
$m$ in the label for $F_j$.  

First observe that the labels in $F_j$ must be weakly increasing at $u$, 
since the latter label
$m$ is the largest value available.  Suppose
there is a maximal face $\sigma\in F_j \cap (\cup_{i<j} F_i)$ which omits 
$u$, and let $\sigma = \hat{0}=v_0 < v_1 < \cdots < v_r < v_{r+1}=\hat{1}$.  
By the maximality of $\sigma $ along with the fact that $F_j$ is increasing
at $u$, the labels on $F_j$ must be 
weakly increasing from $v_i$ to $v_{i+1}$ for
$0\le i\le r$, since otherwise some interval has a descent which 
could be omitted from $F_j$ to obtain a codimension one
face $\tau \in F_j \cap (\cup_{i<j} F_i)$ such that $u\not \in \tau$ and 
$\sigma $ is strictly contained in $\tau $ (contradicting $\sigma $ being
maximal).  Assume that $\sigma $
is maximal in $F_j\cap (\cup_{i<j} F_i)$, $\sigma $ omits $u$, and 
that $F_j$ is weakly increasing between any two elements of $\sigma $.
Let $F_{i'}$ be one of the facets that is lexicographically smalier 
than $F_j$ and contains $\sigma $.  For such an $F_{i'}$ to exist, we 
need there to exist a permutation $\pi $ permuting the row values such
that $F_j |_{supp(\sigma )} = \pi F_{i'}|_{supp(\sigma )} $ and such that
$F_{i'}$ is lexicographically smaller than $\pi F_{i'}$.
This guarantees us the following properties of $\sigma $:

\begin{enumerate}
\item
$\sigma $ must skip one or more intervals of $F_j$ such that two of the 
labels $R_1,R_2\in \{ 1,\dots ,m\} $ each
first appear in the first of these intervals.  Let us assume $R_1<R_2$.
\item
Within each of the intervals of $F_j$ skipped by $\sigma $
the labels $R_1$ and $R_2$ appear an equal (nonzero) number of times
\item
On each of the intervals skipped by $\sigma $, $R_1$ is the 
smallest label and $R_2$ is the largest label 
\item
The first appearance of $R_2$ labeling a covering relation 
$v\prec w$ such that $v,w\in \sigma $ is at a lower rank 
than the first such appearance of $R_1$
\item
$\sigma $ is missing at least one interval below $v$
\end{enumerate}

Observe that a face $\sigma $ 
meeting the above conditions cannot omit $u$ because that would imply 
that $R_2=m$, contradicting the fact that $m$ must later appear as a 
label on the 
covering relation $v\prec w$, since we chose $u$ to immediately 
precede the highest rank appearance of $m$ as a label.
\EOP

\begin{rk}
This argument generalizes immediately to the link of any face which omits
a single interval upon which the largest label appears more than once.  
\end{rk}

\medskip
We show next how to relax this requirement on the largest label to the
requirement that some label appear more than once. 

\begin{prop}\label{intcollapse}
Let $\sigma $ be a face that omits a single interval $T<S$ such that
$S=T\cup S'$ and some letter in $S'$ appears with multiplicity greater
than one.  Then $\rm{lk}(\sigma )$ is collapsible.
\end{prop}

\proof
Notice that any saturated chain $F_j$ for $j>1$ in 
$lk(\sigma )$ has a cone point in 
$F_j \cap (\cup_{i<j} F_i)$ located at
the rank immediately before the last appearance of the 
largest label which does not appear exclusively in a rooted chain
$T \prec u_1 \prec \cdots \prec u_r$ with strictly decreasing labels.
For example, we claim that a saturated chain labeled $43212$ has a cone
point immediately before the second appearance of the label $2$, since the
labels $3$ and $4$ are eliminated by our requirement.  The 
argument of the preceding theorem carries over easily to verify that
this is a cone point, and the existence of such a rank follows from our
requirements on $\sigma $.
\EOP

\begin{qn}
Can collapsibility also be deduced for links which are not a 
single interval, when at least one (or perhaps all) the intervals in the 
link satisfy the conditions of Proposition ~\ref{intcollapse}?  Notice
that links of faces in quotient complexes are not simply joins of links
of faces each omitting a single interval.
\end{qn}

\section{Partitioning $\Delta (B_{lm} )/S_l\wr S_m$}

The complex $\dwr $ is not shellable for $l>2,m>1$, but this section 
provides a partitioning for each pair $l,m$.  This will involve
a labeling that is quite a bit different from the one appearing in
earlier sections.

\begin{defn}
A {\bf partitioning} of a pure boolean cell complex $\Delta $
is an assignment of a face $G_i$ to each facet $F_i$ 
so that the boolean upper intervals $[G_i,F_i]$ partition the set of
faces in $\Delta $, i.e. so that $\Delta $ is
a disjoint union of boolean algebras $[G_1,F_1]\cup\cdots\cup [G_s,F_s]$
whose maximal elements are the facets of $\Delta $.  
\end{defn}

\medskip
A partitioning of a pure, balanced complex $\Delta $ gives a combinatorial
interpretation for the flag $h$-vector, namely each coordinate 
$h_S(\Delta )$ counts minimal faces $G_i$ of support $S$ in the 
partitioning.  We begin with
an example of how to partition a certain link which is not 
Cohen-Macaulay over the integers,
before turning our attention to the entire complex $\dwr $.  Throughout
this section, we use the isomorphism $$\dwr \cong (\Delta(B_{lm})/S_l
\times\cdots\times S_l )/S_m$$ which allows us to view vertices as subsets
of $\{ 1^l,\dots,m^l\} $ modulo an $S_m$-action permuting values.  We refer
to each of the $m$ values as a {\bf row}, motivated by the description of 
$S_l\wr S_m $ following Example 3.1.

\begin{figure}[h]
\begin{picture}(200,120)(-20,0)
\psfig{figure=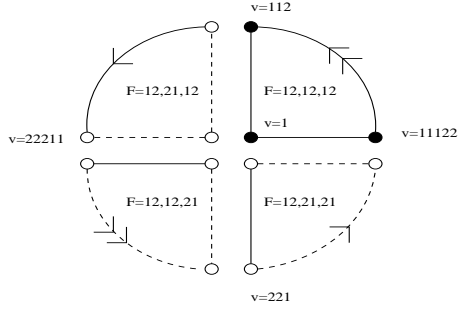,height=4cm,width=6cm}
\end{picture}
\caption{A partitioning for $\real P_2$}
\label{rp2}
\end{figure}

\begin{examp}
{\rm Consider the quotient complex $\Delta (B_6)/S_3 \wr S_2$ and the face
$F=\emptyset < \{ 1,2\} < \{ 1^2,2^2\} < \{ 1^3,2^3\} $.  Notice
that $lk F \cong \real P_2$, as depicted in Figure 
~\ref{rp2} with the usual boundary identifications.  Here, we 
represent the four
facets by $3$-tuples $(\sigma_1,\sigma_2,\sigma_3 )$ of permutations
in $S_2$, written in one-line notation, with the requirement that
$\sigma_1 = 12$.  Vertices in the link are subsets $S$ of $\{ 1^3,2^3\}$ with
the multiplicity of 1 and 2 differing by one, with the identification
$S=(12)S$.
We assign minimal faces to facets as 
follows: $(12,12,12)$ is assigned the empty chain, $(12,12,21)$ is
assigned the chain $\{ 1\} < \{ 1^2, 2^3 \}$, $(12,21,12)$ is assigned
the chain $\{ 1,2^2 \} < \{ 1^3 ,2^2 \}$ and $(12,21,21)$ is assigned
the chain $\{ 1\} < \{ 1,2^2\} $.  In Figure ~\ref{rp2},
vertices and edges of a facet that are assigned
by the partitioning to a different facet are depicted by hollow circles
and dashed edges, respectively.  } 
\end{examp}

\medskip
The partitioning in the above example generalizes to 
$\rm{lk}(
\emptyset < \{ 1,\dots ,m \} < \cdots < \{ 1^l,2^l,\dots ,m^l \} )$ in
$\dwr $ by representing facets by $l$-tuples $(\sigma_1,\dots ,
\sigma_l)\in (S_m)^l$ such that $\sigma_1$ is the identity permutation,
and including in the minimal face associated to 
$(\sigma_1,\dots ,\sigma_l)$ exactly the ranks $im+j$ for $0\le i< l,
0<j<m$ such that $\sigma_i^{-1}(\sigma_{i+1}(j))>\sigma_i^{-1}
(\sigma_{i+1}(j+1))$ (letting $\sigma_0 = \sigma_m$).
We omit the justification of this construction,
instead showing how to partition the entire complex 
$\dwr $ in a related fashion and verifying the validity of that
construction.

The partitioning for $\dwr $ will make use a notion of ascents and 
descents in the facets, based on a labeling for the covering 
relations.
This labeling will give a unique increasing chain on each
interval, and the descents will specify which ranks to include
in the minimal faces assigned to facets.  However, the labeling will not 
give a lexicographic shelling for three reasons: (1) the labeling is not
a chain-labeling, because the 
label assigned to a covering relation depends not only on the chain below, 
but also
on whether the label is being compared with the one below it or above it in
the chain,
(2) the increasing chain is not always lexicographically smallest on
an interval and (3) we define increasing to mean each pair of consecutive
labels is increasing, but because of (1), this 
is not the same as the entire chain increasing.  

We will use a 
permutation $\sigma $ that evolves as we proceed upward from $\hat{0} $ 
to $\hat{1} $ in a saturated chain to play a similar role to the 
$l$-tuple $\sigma_1,\dots ,\sigma_l $ that appeared immediately after 
Example 5.1.  For each vertex in a saturated chain orbit, 
$\sigma $ provides an ordering on the rows 
from which letters are chosen.  Since the choice of permutation
$\sigma $ depends both on the saturated chain orbit and also 
on the rank within that chain, we will denote the 
permutation at rank $i$ by $\sigma_i (C)$ when the rank seems necessary to
clarify meaning, and we will sometimes omit the rank-indicator.

The permutation $\sigma\in S_m$ is initialized to the identity,
and evolves as we proceed from $\hat{0} $ to $\hat{1}$ in a saturated chain
by moving a row $R$ in front of all the rows that are currently similar to
it (as defined below)
whenever a covering relation $T\subset S$ enlarges a set $T$ to $S$ by adding
an element from row $R$.

Before we define row similarity, let us establish a notion of similarity 
block, though its definition will be inductively intertwined with the 
definition of row similarity.

\begin{defn}
A series of consecutive covering relations $u_0\prec\cdots\prec u_{st}$ is
called a {\bf similarity block} if there is some collection of rows
$R_1,\dots ,R_t$ that are similar in $u_0$ and that have each been chosen
the same 
number of times in the saturated chain from $\hat{0}$ to $u_0$ such that
for $0\le i< t$ the covering relations $u_{is}\prec\cdots\prec u_{(i+1)s}$
all insert copies of the row $R_{i+1}$.  
\end{defn}

\medskip
Notice that the second requirement on $u_0$ ensures that similarity blocks
are non-overlapping, and also note that the rows 
$R_1,\dots ,R_t$ may be listed in any order.

\begin{defn}
Let us define {\bf similarity} of rows recursively as follows: all of 
the rows are similar in a saturated chain $C$ at $\hat{0}$.  A collection
of rows $R_1,\dots ,R_t$ which are similar at $u$ will still be similar
at $v$ for $u<v$ if every time any one of the rows $R_i$ 
appears in the interval from $u$ to $v$, 
it appears as part of a similarity block involving the rows $R_1,\dots ,R_t$
(though this similarity block might continue beyond $v$ or begin prior to
$u$).
\end{defn}

\medskip
Thus, fewer and fewer rows will be similar to a fixed row $R$ as we proceed
from $\hat{0} $ to $\hat{1}$.  At the point $u$ when 
rows $R$ and $R'$ cease to 
be similar because of $R$ appearing in a similarity block that does not 
contain $R'$ we have $\sigma_{rk(u)}(R) < \sigma_{rk(u)}(R')$ (and more 
generally we have $\sigma_j(R)<\sigma_j(R')$ for $j\ge \rm{rk}(u)$).

\begin{examp}
Consider the saturated chain orbit which sequentially chooses elements
from three rows in the following order: $112221132333321$.  Notice that
similarity of rows 1 and 2 lasts until the covering relation inserting the 
first 3; row 3 ceased to be similar to the other two rows at the 
covering relation inserting the third 2.  Listing those permutations 
$\sigma_i(C)$ in one-line notation that differ from $\sigma_{i-1}(C)$, we
get $\sigma_{init}(C)=123, \sigma_5=213$ and $\sigma_7=123$.
\end{examp}

The eventual row order $\sigma_{final} $ 
is used to determine descents from wrap-around.
At any particular rank, $\sigma $ reflects the 
partial evolution from the identity permutation based on row insertion
up to this point.  

In analogy to our use of $\sigma_i^{-1}\circ
\sigma_{i+1} $ (in which we let $\sigma_0=\sigma_{final}$) following 
Example 5.1, let us
now consider the renormalized permutations $\rho_r(C) = 
\sigma^{-1}_{final}(C)\circ \sigma_r(C)$.
When a covering relation $T\subseteq S$ adds to $T$
an element from a new row, by convention let us choose this element to
come from
the earliest row not yet chosen.  The label for each insertion is the 
pair $(i,\rho(j))$ where $j$ is the row being chosen
and $i$ is the number of times row 
$j$ has been chosen so far in the chain (including its current selection);
the permutation $\rho $ is evaluated either at $T$ (when comparing to 
a higher covering relation $S\subseteq S'$) or at $U$ (when 
$T\subseteq S$ is being compared to a lower covering relation
$U\subseteq T$).

\begin{defn}
The {\bf relative transpose} order (cf. [He, p. 25]) on labels $(i,\rho (j))$
is a rule for comparing two consecutive covering relation labels in a 
saturated chain.  We compare
covering relations $u\prec v$ and $v\prec w$ by comparing 
their labels $(i_1,\rho_{\rm{rk}(u)}(j_1))$ and 
$(i_2,\rho_{\rm{rk}(u)}(j_2))$, and we say that
$(i_1,\rho_{\rm{rk}(u)}(j_1)) < (i_2,\rho_{\rm{rk}(u)}(j_2))$ 
if $i_1 < i_2$ or if $i_1 = i_2$ and $\rho(j_1) < \rho(j_2)$.
\end{defn}

\medskip
This edge-comparison rule is designed for the sole purpose of
specifying which ranks are
ascents and which are descents.  We call a chain
{\bf increasing} on an interval if it has no descents in the relative
transpose order on that interval, and likewise a {\bf decreasing chain} must
have all descents on the interval.
Our partitioning assigns minimal faces $G_j$ to the 
facets $F_j$ by including in $G_j$ the ranks of the descents in $F_j$ in
the relative transpose order.

\begin{thm}
This assignment of minimal faces to facets gives a partitioning of
$\Delta (B_{lm})/S_l\wr S_m$.
\end{thm}

\proof 
To ensure that our assignment of minimal faces to facets
gives a partitioning, we must check (1) that every 
face belongs to at least one interval $[G_i,F_i]$ and (2)
that no face belongs to multiple intervals.  To verify (1), 
we describe in Proposition 5.1 
how to extend any face $F$ to a facet $F_j$ whose minimal
face $G_j$ is contained in $F$; it suffices to show that $supp(G_j) 
\subseteq supp(F)$ since $F$ and $G_j$ are both faces of $F_j$.  The
second (much easier) claim is confirmed in Proposition 5.2.
\EOP

\begin{prop}\label{cover}
Every face $F$ is contained in an interval $[G_j,F_j]$.
\end{prop}
\proof
Let us describe how to extend each face
$F$ to a facet $F_j$ in such a way that descents in the relative transpose
order on labels of $F_j$ only occur at ranks in the support of $F$.  We 
obtain such an $F_j$ by (1) extending $F$ to a facet $\overline{F}$
in such a way that the extension of 
each interval of $F$ would be increasing (in the relative transpose order)
if $\sigma_{final}(\overline{F})$ 
were the identity permutation, then (2) relabeling 
the rows (since this preserves
the facet orbit) so that the relabeling of $\sigma_{final}(\overline{F})$ 
written in one-line notation is the identity permutation,
then (3) restricting to the resulting representation of the face orbit $F$ 
(which is no longer in standard form), and finally (4) taking $F_j$ to be
the increasing extension 
of this representation of $F$, using the fact (to be confirmed in Lemma
~\ref{finalid}) that $\sigma_{final }(F_j)$ is
the identity permutation.  Example ~\ref{extendrep} provides an example of
this process; notice that $F_j \ne \overline{F} $ in the example, and
that the relabeling of $\overline{F}$ 
has the same set of descents in the relative transpose order as $\overline{F}$
did.
Once we check that $\sigma_{final}(F_j)$ equals the 
identity permutation, we will know that $F_j$ is increasing on 
every interval of $F$, implying $supp(G_j)\subseteq supp(F)$.  
\EOP

\begin{examp}\label{extendrep}
\rm{Let $F= \{ 1^2,2\} < \{ 1^2,2^3\}$ in $\Delta(B_6)/S_3\wr S_2$, 
so then $\overline{F}$ is the saturated
chain in which row elements are inserted in the following order: 
$112221$.  Notice that $\sigma_{final}(\overline{F}) $ is the adjacent
transposition $21$.  Thus, we relabel by swapping 
the values $1$ and $2$, so the 
relabeled representation of $\overline{F} $ is $221112$.  This restricts to
the new representation for $F$ as $F_{relabel} = \{ 2^2,1\} < \{ 2^2,1^3\}$, 
which extends to $F_j$ by inserting rows as follows:
$122112$.  Notice that 
$G_j = \{ 1,2^2 \} < \{ 1^3,2^2\} = F$, since $F_j$ has descents at ranks
3,5, and that $F$ belongs to the interval
$[G_j,F_j]$, as desired.}
\end{examp}

\begin{lem}\label{finalid}
Each facet $F_j$ constructed in Proposition ~\ref{cover} has 
$\sigma_{final}(F_j) $ equalling the identity permutation.
\end{lem}

\proof
We will show that $\sigma_{final}(F_j)$ has no inversion pairs.
Suppose the similarity of rows $r$ and $s$ is broken in $\overline{F}$
on the interval $u<v$ for $u,v$ consecutive elements of the chain $F$.  
Let $F_{relabel}$ denote the expression for
$F$ in which the rows are permuted so that
$\sigma_{final}(\overline{F} )$ is relabeled as the identity permutation.
Let us similarly view $u$ and $v$ in this relabeled form.  Because this
relabeling of $\overline{F} $ sends $\sigma_{final }$ to a permutation
with no inversions, we may conclude that in the relabeled pair $u<v$,
that $\sigma_{rk(v)}(r) < \sigma_{rk(v)}(s)$.  Since the relabeled 
$\overline{F}$ is increasing on the relabeled interval $u<v$, we then know
that $v$ has more copies of $r$ than of $s$, and that 
one of the following properties must hold (letting $v<w$ be the interval of
$F$ immediately following $u<v$) to ensure that there is no similarity block
for $r$ and $s$ beginning on the interval $u<v$ and concluding on the 
interval $v<w$:
\begin{enumerate}
\item
the number of new copies of $r$ in $v<w$ is larger than the number of 
new copies of $s$ on the interval $v<w$
\item
the interval $u<v$ also inserts letters with larger labels than $r,s$, 
implying that these are inserted after the copies of $r$ and $s$, preventing
the continuation of a similarity block to the interval $v<w$
\item
some row $t$ which has smaller value than $r$ or $s$
(and so would precede 
any copies of $r$ or $s$ in the interval $v<w$) is inserted
in the interval $v<w$, again preventing the continuation of a similarity
block to $v<w$
\end{enumerate}

\medskip
One may easily check that these properties carry over to the intervals
$u<v,v<w$ in $F_j$ by virtue of 
(1) $F_j$ containing the relabeled face $F$, (2) $F_j$ increasing on intervals,
(except possibly from wrap-around) and (3) the fact that similarity of $r,s$
cannot be broken earlier in $F_j$, by virtue of 
the same characterization of how similarity
is broken applied to the earlier intervals.  We conclude that the
permutation $\sigma_{final}(F_j)$ has exactly the same inversion pairs as
the relabeling of $\sigma_{final}(\overline{F} )$, so $\sigma_{final}(F_j)$
is the identity permutation.
\EOP

\medskip
It is easy to check that each face is included only once in the partitioning.  

\begin{prop}
There is no overlap among the intervals $[G_j,F_j]$.
\end{prop}
\proof
If $F\in [G_j,F_j]$, then $F_i$ 
must be increasing in the relative-transpose order on
each interval of $F$.  The only possible
flexibility in how to extend $F$ to $F_j$ 
comes from the choice of presentation of $F$ prior to taking its increasing
extension, but at most
one such choice will yield $\sigma_{final} $ which equals the 
identity permutation, as needed to avoid descents from wrap-around.
\EOP

\medskip
As a reality check, we 
computed that the determinant of the incidence matrix for the partitioning of 
$\rm{lk}(\emptyset < \{ 1,2\} < \{ 1^2,2^2\} < \{ 1^3,2^3\} )$ is 2 and 
that the incidence matrix $M$ for the partitioning of 
$\rm{lk}(\emptyset < \{ 1,2\} < \{ 1^2,2^2\} < \{ 1^3,2^3\} 
< \{ 1^4,2^4 \} )$ has $det(M)=8$, 
consistent with the fact that $\real P_n$ only has local 2-torsion.

\begin{qn}
Is the incidence matrix $M$ for this partitioning
nonsingular over $\integ /p\integ $ for all $p>2$?  
If so, then 
the partitioning would give a basic set for the subring 
$k[x_1,\dots ,x_{lm}]^{S_l\wr S_m}$ of polynomials that are invariant
under the action of $S_l\wr S_m$ for $char(k)>2$, 
by results of [GS] about transferring basic sets.  This would imply 
Cohen-Macaulayness for $char(k)>2$.  
\end{qn}

\medskip
We suspect that this question has a negative answer.
Notice that $M$ is nonsingular over $\integ /p\integ $ if and only if
$p$ does not divide the determinant of $M$.  
The incidence matrix $M$ for our 
partitioning for $\rm{lk}(\emptyset < \{ 1,2,3\} <
\{ 1^2,2^2,3^2 \} < \{ 1^3,2^3,3^3\} )$ satisfies
$det(M)=2^3\cdot 3^5$, and so is singular over $\integ /3\integ $, suggesting
the distinct possibility of local $3$-torsion.

\section{Appendix (by Vic Reiner)}
We wish to prove Theorem~\ref{appendix-result}.  For this purpose, we
introduce some notation, which mostly follows that of [GS]:
$$
\begin{aligned}
R:=&\text{ Stanley-Reisner ring for the Boolean algebra }B_n - \{\emptyset\} \\
 =&k[\,\, y_S: \emptyset \ne S \subseteq [n]\,\,]/ I \\
 & \text{ where }I\text{ is the ideal generated by all products }y_S y_T \\
 & \text{ with }S,T\text{ incomparable subsets of }[n]\\
R':=& k[x_1,...,x_n]
\end{aligned}$$
$$
\begin{aligned}
G =& \text{ a subgroup of the symmetric group }S_n, \\
   &\text{ acting on both $R, R'$ by permuting subscripts.}\\
R^G, (R')^G &\text{ the corresponding invariant subrings.}\\
T:=& \text{ the transfer map }R \rightarrow R'\text{ from [GS]}, 
                    \text{ mapping  }y_S \mapsto \prod_{i \in S} x_i,\\
                   & \text{ then extending multiplicatively to non-vanishing monomials}\\
                   & \text{ in $R'$, then further extending 
                    $k$-linearly to all of }R'.\\
\theta_i : = &\sum_{S:|S|=i} y_S  \in R\\
k[\theta ] : =&k[\theta_1,...,\theta_n] \subset R\\
e_i :=&\text{ the }i^{th}\text{ elementary symmetric function in }x_1,\ldots,x_n \\
     =&T(\theta_i)\\
k[e]:=& k[e_1,...,e_n] \subset R'\\
     =& T(k[\theta]).
\end{aligned}
$$
 
\begin{thm}
If $R^G$ is a Cohen-Macaulay ring, then $(R')^G$ is also a Cohen-Macaulay
ring.  
\end{thm}

\proof
If $R^G$ is Cohen-Macaulay, then the h.s.o.p. \!\! $\theta_1,\dots ,\theta_n$
is a regular sequence, so $R^G$ is a free module
over the polynomial ring $k[\theta]:=k[\theta_1,...,\theta_n]$.
Furthermore, we can choose a basis for this free module consisting of
elements $\eta_1,\dots ,\eta_t$ which are homogeneous with respect to the
fine $\natnum^n$-grading on $R^G$ (choosing any $\eta_i's$ which are
finely homogeneous liftings of a $k$-vector space basis for
$R^G / ( k[\theta]_+ ) $ will work).

We wish to show that $T(\eta_1),\dots ,T(\eta_t)$ comprise a $k[e]$-basis for 
$(R')^G$ as a free $k[e]$-module, which would then show that $(R')^G$ is
Cohen-Macaulay.
We first argue by a comparison of Hilbert series
that one only needs to show that $T(\eta_1),\dots ,T(\eta_t)$
span.  Since $T$ is a $G$-equivariant $k$-vector space isomorphism (but not
a ring isomorphism!) from $R$ to $R'$, it restricts to a $k$-vector space
isomorphism from $R^G$ to $(R')^G$.  If, for the moment, we 
coarsely $\natnum$-grade $R^G$
by applying the usual specialization to its fine $\natnum^n$-grading
(i.e. so that $y_S$ has degree $|S|$), then
$T$ also respects the polynomial gradings on each side.  This
implies $R^G$ and $(R')^G$ have the same Hilbert series.  
Hence the fact that $\eta_1,\ldots,\eta_t$ form a free $k[\theta]$-basis
for $R^G$ implies that the degrees of $T(\eta_1), \dots ,T(\eta_t)$
are such that there are the right number of
$k[e]$-linear combinations of them in each degree to form a basis
of $(R')^G$.  If we can show that $T(\eta_1), \dots ,T(\eta_t)$ do span $(R')^G$
as a $k[e]$-module, we would then know that these $k[e]$-linear combinations
give a $k$-basis in each degree, so they would form a $k[e]$-basis for $(R')^G$.
 
For the spanning argument, since $(R')^G$ is spanned as a $k$-vector space
by $G$-orbit sums $G(x^{\alpha})$ of monomials $x^{\alpha} \in R'$, 
we only need
to show that such elements are in the $k[e]$-span of the $T(\eta_i)'s$.
Let $G(x^{\alpha })$ be any such $G$-orbit sum.
Let $T^{-1}(G(x^{\alpha}))$ have an expression in $R^G$ as follows:
 
$$ T^{-1}(G(x^{\alpha})) = \sum_i \eta_i p_i(\theta) \eqno(*)$$
for some polynomials $p_i$ in the $\theta's$.
 
We will show that
 
$$    G(x^{\alpha }) - \sum_i T(\eta_i) p_i(e) \eqno(**)$$
is a sum of monomials $x^{\beta }$ whose ``shapes''
(as defined in [GS, p.178]) are all lower in the dominance order
than the shape of $x^{\alpha }$, using [GS, Lemma 9.1], and then be
done by induction on the dominance order.
 
To see this, note that the shape of $x^{\alpha }$ (and every other monomial
occurring in $G(x^{\alpha})$) is the same as the fine grading of the
element $T^{-1}(G(x^{\alpha }))$, so that in expression (*), we may
assume that every term in the sum has this same $\natnum^n$-grading
(by $\natnum^n$-gradedness of $R^G$).  Then [GS, Lemma 9.1] tells us 
that every monomial one obtains
by multiplying out the terms in the sum in (**) will have shape
less than or equal to that of $x^{\alpha}$ in dominance order, 
and that those whose shapes
match those of $x^{\alpha }$ exactly correspond to the terms in (*), so 
they all cancel with terms in $G(x^{\alpha })$
due to the equality (*).  The shapes of the remaining
non-cancelling monomials in (**) are 
all strictly lower in dominance order.
\EOP

\section*{Acknowledgments}

The author thanks Vic Reiner for telling her about the relationship
between quotient complexes and subrings of invariant polynomials, and
also for suggesting $\Delta (B_{2m})/S_2\wr S_m$ as a candidate for her
lexicographic shellability criterion, since this also led to
her study of $\dwr $.

\end{document}